\documentclass[12pt]{amsart}
\input diagrams
\diagramstyle[PostScript=dvips]

\newcommand{\cal}{\mathcal}
\renewcommand{\Bbb}{\mathbb}

\newcommand{\ch}{\operatorname{ch}}
\newcommand{\td}{\operatorname{td}}

\newcommand{\DD}{{\cal D}}
\newcommand{\KK}{{\cal K}}
\newcommand{\MM}{{\cal M}}

\newcommand{\mgn}{\overline{\MM}_{g,n}}
\newcommand{\mgnr}{\mgn^{1/r}}

\newcommand{\cv}{c^{1/r}}            

\newcommand{\Cone}{\operatorname{Cone}}

\newcommand{\TT}{{\cal T}}

\newcommand{\lan}{\langle}
\newcommand{\ran}{\rangle}

\newcommand{\CC}{{\cal C}}

\newcommand{\Tr}{\operatorname{Tr}}

\newcommand{\si}{\sigma}

\newcommand{\ga}{\gamma}
\newcommand{\de}{\delta}

\newcommand{\We}{\Lambda}
\renewcommand{\ker}{\operatorname{ker}}
\newcommand{\im}{\operatorname{im}}

\numberwithin{equation}{section}

\newtheorem{thm}{Theorem}[section]

\newtheorem{lem}[thm]{Lemma}

\newenvironment{rem}{\vspace{3mm}\noindent
{\bf Remark.}}{\vspace{3mm}}

\theoremstyle{definition}

\newcommand{\Pf}{\noindent {\it Proof}}
\newcommand{\id}{\operatorname{id}}

\newcommand{\ra}{\rightarrow}

\renewcommand{\AA}{{\cal A}}
\newcommand{\FF}{{\cal F}}

\newcommand{\LL}{{\cal L}}

\newcommand{\OO}{{\cal O}}

\newcommand{\Hom}{\operatorname{Hom}}
\newcommand{\Ext}{\operatorname{Ext}}

\newcommand{\Res}{\operatorname{Res}}

\renewcommand{\a}{\alpha}
\renewcommand{\b}{\beta}
\newcommand{\om}{\omega}
\newcommand{\De}{\Delta}
\newcommand{\la}{\lambda}

\newcommand{\Z}{{\Bbb Z}}
\newcommand{\Q}{{\Bbb Q}}

\newcommand{\wt}{\widetilde}

\newcommand{\sub}{\subset}
\newcommand{\ed}{\qed\vspace{3mm}}

\newcommand{\II}{{\cal I}}
\newcommand{\ev}{\operatorname{ev}}

\begin{document}

\title{Witten's top Chern class on the moduli space of higher
spin curves}

\author{Alexander Polishchuk}
\address
{Department of Mathematics, Boston
University, Boston, MA 02215} 
\email{apolish@math.bu.edu}
\thanks{This work was partially supported by NSF grant
DMS-0070967.}

\begin{abstract} 
We prove that the algebraic Witten's ``top Chern class" constructed
in \cite{PV} satisfies the axioms for the spin virtual class formulated
in \cite{JKV}.
\end{abstract}
\maketitle

This paper is a sequel to \cite{PV}.
Its goal is to verify that the 
{\em virtual top Chern class\/} $\cv$ in the Chow group of the moduli
space of higher spin curves $\mgnr$, constructed in \cite{PV},
satisfies all the axioms of {\em spin virtual class\/} formulated 
in~\cite{JKV}. Hence, according to~\cite{JKV},
it gives rise to a cohomological field theory
in the sense of Kontsevich-Manin~\cite{KM}. 
As was observed in \cite{PV}, the only non-trivial axioms that have
to be checked for the class $\cv$ are two axioms that we call
{\it Vanishing axiom} and {\it Ramond factorization axiom}.
The first of them requires $\cv$ to vanish on all the
components of the moduli space $\mgnr$, where one of the markings is
equal to $r-1$. The second demands vanishing of the push-forward of
$\cv$ restricted to the components of the moduli space 
corresponding to the so called Ramond sector, under some natural
finite maps. 

Recall that the virtual top Chern class is a crucial ingredient
in the generalized Witten's conjecture formulated in~\cite{W1}, \cite{W2}.
The original index-theoretic construction of this class sketched
by Witten was recently extended to the compactified moduli
space by T. Mochizuki \cite{Mo} who also showed that the obtained class
satisfies the axioms of
\cite{JKV}. The algebraic construction of \cite{PV} gives a class
in the Chow group with rational coefficients (and
axioms are satisfied on the level of Chow groups).
Presumably, the algebraic construction induces the same class in cohomology of
$\mgnr$ as the analytic construction.

It is interesting to note that the class $\cv$ is constructed as
a characteristic class of certain supercommutative $\Z/2\Z$-graded dg-algebra
over $\mgnr$ equipped with an odd closed section (where the entire data is defined up to
quasi-isomorphism). This resembles Kontsevich's approach to the construction of the 
virtual fundamental class (see \cite{Kon}). One may hope that both constructions can be embedded
into a more general framework involving dg-spaces. This would be in agreement with the
philosophy of derived moduli spaces promoted in \cite{Kon}, \cite{CK1}, \cite{CK2}.
 
The paper is organized as follows. In section \ref{homsec} we prove two
identities for localized Chern characters of specific $\Z/2\Z$-graded complexes.
In section \ref{vanishsec} we deduce Vanishing axiom from the first identity and
in section \ref{Ramondsec} we deduce Ramond factorization axiom from the second identity.

\medskip


\section{Some $\Z/2\Z$-graded homological algebra}\label{homsec}

Recall (see \cite{PV})
that for every $\Z/2\Z$-graded complex $(V^{\bullet}=V^+\oplus V^-,d)$
of vector bundles on a scheme $X$,
which is strictly exact off a closed subset $Z\subset X$, the
graph-construction associates the {\it localized Chern character}
$\ch_Z^X(V^{\bullet})\in A^*(Z\ra X)_{\Q},$
where $A^*(Z\ra X)_{\Q}$ is the bivariant Chow group with rational coefficients.
(the original construction given in \cite{BFM} or \cite[Ch.\ 18]{F}
deals with $\Z$-graded complexes). Here we use the following
terminology from \cite{PV}:
$(V^{\bullet},d)$ is {\it strictly exact} if it is exact and
$\im(d)$ is a subbundle of $V^{\bullet}$. Note that if a $\Z/2\Z$-graded
complex is homotopic to zero then it is strictly exact (since in this
case $\im(d)$ is a direct summand of $V^{\bullet}$).
Witten's top Chern class is constructed in \cite{PV}
by slightly modifying the localized Chern class of the complex corresponding
to the action of an isotropic section of an orthogonal bundle
on a spinor bundle, where the relevant orthogonal data is constructed
using the higher spin structure on the universal curve over $\mgnr$
(we will recall this construction in section \ref{vanishsec}).

We are going to prove two identities for the localized Chern character
that will be main ingredients for the proof of Vanishing axiom
and Ramond factorization axiom respectively.
In both cases the identities hold on the level of $K$-theory (of complexes
that are strictly exact off a closed subset).

\begin{lem}\label{mainlem}
Let $d(\la):V^{\bullet}\ra V^{\bullet}[\la]$ be an odd
endomorphism of a $\Z/2\Z$-graded bundle $V^{\bullet}$ on $X$ of the form
$d(\la)=d_0+d_1\la+\ldots+d_{r-1}\la^{r-1}$ for some $r\ge 2$,
depending on a formal parameter $\la$ (commuting with everything).
Assume that $d(\la)^2=\la^r$. Then
$\ch_Z^X(V^{\bullet},d_0)=0$ for every closed subset $Z\sub X$ such that
$(V^{\bullet},d_0)$ is strictly exact off $Z$.
\end{lem}

\Pf . We can extend $d(\la)$ to a $\OO_X[\la]$-linear endomorphism
of $V^{\bullet}[\la]$. Let $W^{\bullet}=V^{\bullet}[\la]/(\la^r)$
with the odd endomorphism $d_W:W^{\bullet}\ra W^{\bullet}$
induced by $d(\la)$. Then $d_W^2=0$ and there is a natural
$r$-step filtration on the complex $(W^{\bullet},d_W)$ with all
consequtive quotient-complexes isomorphic to $(V^{\bullet},d_0)$.
Thus, it is enough to prove that $(W^{\bullet},d_W)$ is strictly
exact (everywhere). We claim that in fact this complex is homotopic
to zero. Indeed, for every interval of integers $[a,b]$ let us denote
$V^{\bullet}[\la]_{[a,b]}=\oplus_{i=a}^b V^{\bullet}\la^i$.
Let us denote by $d'$ and $d''$ the following components of the
restriction of $d(\la)$ to $V^{\bullet}[\la]_{[0,r-1]}$:
$$d(\la):V^{\bullet}[\la]_{[0,r-1]}\rTo^{(d',d'')}
V^{\bullet}[\la]_{[0,r-1]}\oplus V^{\bullet}[\la]_{[r,2r-1]}.$$
Note that $d'=d_W$ upon the natural identification of
$V^{\bullet}[\la]_{[0,r-1]}$ with $W^{\bullet}$. Also, the
image of $d''$ is contained in $V^{\bullet}[\la]_{[r,2r-2]}$.
Extending $d'$ and $d''$ to $\OO_X[\la^r]$-linear endomorphisms of
$V^{\bullet}[\la]$ we can write $d=d'+d''$.
Then the
condition $d(\la)^2=\la^r$ implies that $d'd''+d''d'=\la^r$ on
$V^{\bullet}[\la]_{[0,r-1]}$.
Hence
$$h=d''/\la^r:V^{\bullet}[\la]_{[0,r-1]}\ra V^{\bullet}[\la]_{[0,r-1]}$$
gives a homotopy between the identity and zero endomorphisms of
the complex $(V^{\bullet}[\la]_{[0,r-1]},d')\simeq (W^{\bullet},d_W)$.
\ed

\begin{rem} The above lemma admits the following generalization:
if the differential $d(\la)=d_0+d_1\la+\ldots+d_{r-1}\la^{r-1}$
as above satisfies $d(\la)^2=f(\la)$ for some polynomial $f$ of
degree $r$ then
$$\sum_{z: f(z)=0}m_z\cdot\ch_Z^X(V^{\bullet},d(z))=0$$ 
where $m_z$ is the multiplicity of a root $z$,
$Z\subset X$ is a closed subset such that all the complexes
$(V^{\bullet},d(z))$ (for $f(z)=0$) are strictly exact off $Z$.
\end{rem}

\begin{lem}\label{mainlem2}
Let $d:V^{\bullet}\ra V^{\bullet}$ be an odd
endomorphism of a $\Z/2\Z$-graded bundle $V^{\bullet}$ on $X$ 
such that $d^2=-(f_1\ldots f_r)\cdot\id_{V^{\bullet}}$, 
where $f_1,\ldots,f_r$ are functions on $X$. For every
$i=1,\ldots,r$ let us introduce the differential
$d_i$ on $V^{\bullet}\oplus V^{\bullet}[1]$ by the formula
$$d_i(x,x')=
(d(x)+(\prod_{j\neq i}f_j)\cdot x', -d(x')+f_i\cdot x)$$
where $x\in V^{\bullet}$, $x'\in V^{\bullet}[1]$. Then
$$\sum_{i=1}^r\ch_Z^X(V^{\bullet}\oplus V^{\bullet}[1],d_i)=0$$ 
for every closed subset $Z\sub X$ such that all
$(V^{\bullet}\oplus V^{\bullet}[1],d_i)$ are strictly exact off $Z$.
\end{lem}

\Pf . Let us introduce the differential $D$ on
$W^{\bullet}:=(V^{\bullet}\oplus V^{\bullet}[1])^{\oplus r}$ by the formula
$$D(x_i,x'_i)_{i=1,\ldots r}=(y_i,y'_i)_{i=1,\ldots,r}$$
where $x_i,y_i\in V^{\bullet}$, $x'_i,y'_i\in V^{\bullet}[1]$,
\begin{align*}
&y_i=dx_i+f_{i+1}f_{i+2}\ldots f_r\cdot[x'_1+f_1x'_2+f_1f_2x'_3+\ldots+
f_1\ldots f_{i-1}x'_i]\ \text{for}\ i<r,\\
&y_r=dx_r+x'_1+f_1x'_2+f_1f_2x'_3+\ldots+f_1\ldots f_{r-1}x'_r,\\
&y'_i=-dx'_i+f_ix_i-x_{i-1}\ \text{for}\ i\ge 2,\\
&y'_1=-dx'_1+f_1x_1.
\end{align*}
One can easily check that $D^2=0$.
There is a natural decreasing filtration of $W^{\bullet}$ by
subcomplexes $W^{\bullet}=F^1W^{\bullet}\supset\ldots\supset F^rW^{\bullet}\supset
F^{r+1}W^{\bullet}=0$, where 
$$F^jW^{\bullet}=\{(x_i,x'_i)_{i=1,\ldots,r}: x_1=\ldots=x_{j-1}=0,
x'_1=\ldots=x'_{j-1}=0\}.$$
The associated graded quotients are
$$F^jW^{\bullet}/F^{j+1}W^{\bullet}\simeq (V^{\bullet}\oplus V^{\bullet}[1],d_j),$$
$j=1,\ldots,r$. Therefore,
$$\sum_{i=1}^r\ch_Z^X(V^{\bullet}\oplus V^{\bullet}[1],d_i)=\ch_Z^X(W^{\bullet},D).$$ 
It remains to prove that the complex $(W^{\bullet}, D)$ is strictly exact
on $X$. For this we construct a homotopy $h$ 
between the identity and zero endomorphisms
of $W^{\bullet}$. Namely, we set
$$h(x_i,x'_i)_{i=1,\ldots r}=(y_i,y'_i)_{i=1,\ldots,r},$$
where 
\begin{align*}
&y_i=-[x'_{i+1}+f_{i+1}x'_{i+2}+f_{i+1}f_{i+2}x'_{i+2}+\ldots+
f_{i+1}\ldots f_{r-1}x'_r]\ \text{for}\ i<r-1,\\
&y_{r-1}=-x'_r,\ y_r=0,\\
&y'_1=x_r,\ y'_i=0\ \text{for}\ i\ge 2.
\end{align*}
It is easy to check that $Dh+hD=\id_{W^{\bullet}}$.
\ed

\section{Vanishing axiom}\label{vanishsec}

Henceforward all our schemes are assumed to be quasiprojective over a field $k$.
We assume that $\operatorname{char} k>r$ and that $k$ contains all $r$-th roots of unity.
 
Let $\pi:\CC\ra X$ be a family of prestable curves over a
scheme $X$, and let $\TT$ be a family of rank-one
torsion-free sheaves on $\CC$ equipped with a non-zero
homomorphism $b:\TT^r\ra\om_{\CC/X}$, where $\om_{\CC/X}$ is
the dualizing sheaf of $\pi$. In this situation we defined in
section 5.1 of \cite{PV} the class $c(\TT,b)\in A^{-\chi}(X)_{\Q}$,
where $\chi$ is the Euler-Poincar\'e characteristic of members of the family $\TT$.
To construct this class we consider the map $\tau:S^rR\pi_*\TT\ra\OO_X[-1]$ induced by
$b$ and by the trace map $\Tr:R\pi_*\om_{\CC/X}\ra\OO_X[-1]$. As was proved
in Proposition 4.7 of \cite{PV} there exists a complex $C_0\ra C_1$
of vector bundles on $X$ representing $R\pi_*\TT$ such that the map $\tau$
is represented by the chain map of complexes $S^r[C_0\ra C_1]\ra\OO_X[-1]$.
This chain map corresponds to a morphism of vector bundles $\nu:S^{r-1}C_0\ra C_1^{\vee}$.
We can consider the differential $d:C_0\ra C_1$ and the map $\nu$ as sections of
the pull-backs of $C_1$ and $C_1^{\vee}$ to the total space of $C_0$. Then
$s=(d,\nu)$ will be an isotropic section of the orthogonal vector bundle
$p^*C_1\oplus p^*C_1^{\vee}$ on $C_0$, where $p:C_0\ra X$ is
the projection. Moreover, $s$ vanishes exactly on $X$ embedded into $C_0$ by the zero section.
Then we consider the action of $s$ on the spinor bundle 
$\We^*p^*C_1^{\vee}$. The obtained $\Z/2\Z$-graded complex $(\We^*p^*C_1^{\vee},s)$
is exact outside $X\sub C_0$. Therefore, the localized Chern character of this complex
is an element of the bivariant Chow group $A^*(X\ra C_0)_{\Q}\simeq A^*(X)_{\Q}$. 
The class $c(\TT,b)$
is obtained by multiplying this localized Chern character with the Todd class of $C_1$.
Theorem 4.3 of \cite{PV} assures that $c(\TT,b)$ does not depend on the choices made.

\begin{rem} We can consider 
$(\AA:=p_*(\We^*p^*C_1^{\vee}),\de)$ as a sheaf of $\Z/2\Z$-graded
dg-algebras over $X$, where the differential $\de$ is induced by $d:C_0\ra C_1$.
The action of the isotropic section $s$ on $\AA$ has form $d+\epsilon(e)$,
where $\epsilon(e)$ is the operator of multiplication with the $\de$-closed odd
section $e\in\AA$ corresponding to $\nu:S^{r-1}C_0\ra C_1^{\vee}$.
The proof of Theorem 4.3 in \cite{PV} can be converted into the proof of the fact
that the quasi-isomorphism class of the data $(\AA,e)$
is uniquely determined by $(\TT,b)$.
\end{rem}

Let $\si:X\ra\CC$ be a section
of $\pi$ such that $\pi$ is smooth near $\si(X)$ (a marked point).
By abuse of notation we will denote
by $F\mapsto F(\si)$ the operation of tensoring with the line bundle
$\OO_{\CC}(\si(X))$. The following theorem immediately implies
Vanishing axiom for $\cv$ (Axiom 4 of \cite{JKV}).

\begin{thm}\label{axiom1thm} 
Assume that $b$ factors as a composition
$$\TT^r\stackrel{b_0}{\ra}\om_{\CC/X}(-(r-1)\si)\ra\om_{\CC/X},$$
where $\si^*b_0$ is an isomorphism. Then $c(\TT,b)=0$.
\end{thm}

Let us set $L=\TT(\si)|_{\si}$. The map $\si^*b_0$ gives an isomorphism
$$L^r\wt{\ra}\om_{\CC/X}(\si)|_\si\simeq\OO_X.$$
Since we are working with rational coefficients,
we can replace $X$ by its finite \'etale covering
over which $L$ is trivial (right now we do not need to choose 
a specific trivialization of $L$). 

Let $\pi^{(r)}:\CC^{(r)}\ra X$ denote the relative $r$-th symmetric power
of $\CC$ over $X$. We denote by $\si^r\in\CC^{(r)}$ the $X$-point
corresponding to the relative divisor $r\si(X)$ and by
$\II_{\si^r}\sub\OO_{\CC^{(r)}}$ the ideal sheaf of $\si^r(X)\sub\CC^{(r)}$.
For every coherent sheaf $\FF$ on $\CC$, let 
$\FF^{(r)}$ denote the $r$-th symmetric power of $\FF$, which is a sheaf
on $\CC^{(r)}$. We claim that $b_0$ induces a morphism
\begin{equation}\label{idealsheafmap}
R\pi^{(r)}_*(\II_{\si^r}\otimes(\TT(\si))^{(r)})\ra R\pi_*(\om_{\CC/X}).
\end{equation}
Indeed, let $\De:\CC\ra\CC^{(r)}$ be the diagonal map. Then 
we have a natural morphism
$$\II_{\si^r}\otimes(\TT(\si))^{(r)}\ra
\De_*\De^*(\II_{\si^r}\otimes(\TT(\si))^{(r)})\ra\De_*\TT^r((r-1)\si)\ra
\De_*\om_{\CC/X},$$
where the last arrow is induced by $b_0$. Now the morphism
(\ref{idealsheafmap}) is obtained by applying the functor $R\pi^{(r)}_*$.
Composing (\ref{idealsheafmap}) with the trace map 
$R\pi_*(\om_{\CC/X})\ra\OO_X[-1]$, we get a morphism
$$\wt{\tau}:R\pi^{(r)}_*(\II_{\si^r}\otimes(\TT(\si))^{(r)})\ra \OO_X[-1]$$
that will play a major role in the proof of Theorem \ref{axiom1thm}.
Note that we also have a natural map
$$\iota:R\pi^{(r)}_*(\TT^{(r)})\ra R\pi^{(r)}_*(\II_{\si^r}\otimes(\TT(\si))^{(r)})$$
and an isomorphism $S^r R\pi_*\TT\simeq R\pi^{(r)}_*(\TT^{(r)})$,
such that $\tau=\wt{\tau}\circ\iota$ can be identified with the map
$S^r R\pi_*\TT\ra\OO_X[-1]$ used in the definition of $c(\TT,b)$.

Note that there is an exact triangle
\begin{equation}\label{symextri}
\OO_X[-1]\stackrel{\de}{\ra}
R\pi^{(r)}_*(\II_{\si^r}\otimes(\TT(\si))^{(r)})\ra R\pi^{(r)}_*(\TT(\si))^{(r)}\ra \OO_X
\end{equation}
where we use the canonical trivialization of $L^r$.

\begin{lem}\label{composlem} 
The composition $\wt{\tau}\circ\de$ is the identity map.
\end{lem}

\Pf . This follows immediately from the existence of a natural morphism
of exact triangles
$$
\begin{diagram}
\De^*\II_{\si^r}\otimes\TT^r(r\si) &\rTo & \TT^r(r\si) &\rTo & \De^* (\si^r)_*\OO_X &\rTo\ldots\\
\dTo & & \dTo & & \dTo &\\
\om_{\CC/X} &\rTo &\om_{\CC/X}(\si) &\rTo & \si_*\OO_X &\rTo\ldots
\end{diagram}
$$
and from the fact that the composition 
$$\OO_X[-1]\simeq R\pi_*\si_*\OO_X[-1]\ra R\pi_*\om_{\CC/X}\ra\OO_X[-1]$$
is the identity map.
\ed

We want to realize the maps $\wt{\tau}$ and $\iota$ on the level of complexes in
a compatible way.
We start by realizing the canonical distinguished triangle in 
$D^b(X)$:
$$R\pi_*\TT\stackrel{\a}{\ra} R\pi_*(\TT(\si))\stackrel{\b}{\ra} L
\stackrel{\ga}{\ra} R\pi_*\TT[1]$$
by an exact triple of complexes of vector bundles on $X$.

\begin{lem}\label{extensionlem} 
Let $[C_0\stackrel{d}{\ra} C_1]$ be a complex of vector bundles 
(concentrated in degrees $[0,1]$) representing $R\pi_*\TT$.
Then there exists an extension of vector bundles
\begin{equation}\label{mainexseq}
0\ra C_0\ra\wt{C}_0\ra L\ra 0,
\end{equation}
and a morphism $\wt{d}:\wt{C}_0\ra C_1$ extending $d$, such that
the morphism $\ga: L\ra R\pi_*\TT[1]$
is represented by the chain map of complexes
\begin{equation}\label{chainext}
[C_0\ra\wt{C}_0]\stackrel{(\id,\wt{d})}{\ra} [C_0\ra C_1],
\end{equation}
hence, the complex $[\wt{C}_0\stackrel{\wt{d}}\ra C_1]$
represents $R\pi_*(\TT(\si))$ and the morphisms $\a$ and $\b$
are represented by the natural chain maps
$$[C_0\ra C_1]\ra [\wt{C}_0\ra C_1]\ra L;$$
\end{lem}

\Pf .
Applying the second arrow in
the exact sequence
\begin{equation}\label{exseq}
\Hom(L,C_1)\ra\Hom(L,R\pi_*\TT[1])\ra\Ext^1(L,C_0)\ra\Ext^1(L,C_1)
\end{equation}
to the element $\ga$ we get an extension class 
$e\in\Ext^1(L,C_0)$ which becomes trivial in $\Ext^1(L,C_1)$.
Let 
$$0\ra C_0\ra\wt{C}_0\ra L\ra 0$$
be an extension with the class $e$, $\wt{d}:\wt{C}_0\ra C_1$
be a splitting of its push-out by $d:C_0\ra C_1$.
The element in $\Hom(L,R\pi_*\TT[1])$
represented by the chain map (\ref{chainext}) induces the
same class $e$ in $\Ext^1(L,C_0)$. 
Now the sequence (\ref{exseq})
shows that after changing a splitting $\wt{d}$ by an appropriate
element of $\Hom(L,C_1)$
the chain map (\ref{chainext}) will represent $\ga$.
\ed

Let $C=[C_0\ra C_1]$ be a complex representing $R\pi_*\TT$ and
let $\wt{C}=[\wt{C}_0\ra C_1]$ be the complex representing $R\pi_*(\TT(\si))$
obtained by applying the above lemma. Then the complex
$S^r C$ (resp. $S^r \wt{C}$) represents 
$S^r R\pi_*\TT\simeq R\pi^{(r)}_*\TT^{(\si)}$
(resp. $R\pi^{(r)}_*(\TT(\si))^{(r)}$) and we have a natural surjective map of complexes
$S^r\wt{C}\ra L^r\simeq\OO_X$ induced by the map $\wt{C}_0\ra L$.
Then the kernel complex $\ker(S^r\wt{C}\ra \OO_X)$ represents
$R\pi^{(r)}_*(\II_{\si^r}\otimes(\TT(\si))^{(r)})$ in a way compatible with the exact triangle
(\ref{symextri}).
Moreover, the map $\iota$ is represented by the natural chain map
$S^r C\ra \ker(S^r\wt{C}\ra \OO_X)$.
It remains to choose our data in such a way that $\wt{\tau}$ would be represented by
a chain map $\wt{\tau}:\ker(S^r\wt{C}\ra \OO_X)\ra \OO_X[-1]$.
For this we use the following lemma analogous to Proposition 4.7 from
\cite{PV}.

\begin{lem}\label{complexlem} 
There exists a complex of vector bundles
$C_0\ra C_1$ representing $R\pi_*\TT$, such that 
one has
$$\Hom_{K^b(X)}(E,\OO_X[n])\simeq\Hom_{D^b(X)}(E,\OO_X[n])$$
for $n\le 0$ and $E=\ker (S^r[\wt{C}_0\ra C_1]\ra L)$.
\end{lem}

\Pf . We start with an arbitrary complex of vector bundles
$C'_0\ra C'_1$ representing $R\pi_*\TT$ and then replace
it by the quasiisomorphic complex $C_0\ra C_1$, where
$C_1=\OO_X(-m)^{\oplus N}\ra C'_1$ is a surjection (see \cite{PV}, Lemma 4.6), 
$\OO_X(1)$ is an ample line bundle on $X$
$m$ is an integer (later we will need to choose $m$ sufficiently large). 
The spectral sequence computing
$\Hom_{D^b(X)}(E,\OO_X[n])$ shows that to prove (ii) it suffices to
check the vanishing 
$$H^i(S^j\wt{C}_0^{\vee}(m'))=0$$
for $i>0$, $j<r$, $m'\ge m$.
Since $\wt{C}_0$ is an extension of the trivial bundle by
$C_0$, this would follow from the vanishing
of $H^i(S^jC_0^{\vee}(m'))$
under the same conditions on $i,j,m'$.
We know that for sufficiently large $m$ one has
$$H^{>0}(S^j(C'_0)^{\vee}\otimes
S^{j_1}(C'_1)^{\vee}\otimes\ldots\otimes S^{j_k}(C'_1)^{\vee}(m'))=0$$
for $j+j_1+\ldots+j_k<r$ and $m'\ge m$.
As was shown in Proposition 4.7 of \cite{PV},
this implies that $H^{>0}(S^jC_0^{\vee}(m'))=0$ for $j<r$, $m'\ge m$.
\ed


\noindent
{\it Proof of Theorem \ref{axiom1thm}}.
Let us choose the data $(C_0,C_1,\wt{C}_0,d,\wt{d})$ as in Lemmas
\ref{complexlem} and \ref{extensionlem}. Let us set 
$K:=\ker(S^r[\wt{C}_0\ra C_1]\ra \OO_X)$.
Then the morphism $\wt{\tau}$ is represented by the chain map 
$K\ra \OO_X[-1]$ that corresponds to a morphism
$$\wt{\tau}:S^{r-1}\wt{C}_0\otimes C_1\ra\OO_X$$
such that the composition 
\begin{equation}\label{zerocomp}
\ker(S^r\wt{C}_0\ra\OO_X) \ra S^{r-1}\wt{C}_0\otimes C_1\ra\OO_X
\end{equation}
is zero.

Let $X'\ra X$ be the affine bundle classifying splittings
of the exact sequence (\ref{mainexseq}). Since
the pull-back induces an isomorphism of Chow groups of $X$ and $X'$
we can make a base change of our data by the morphism $X'\ra X$.
Thus, we can assume that the extension (\ref{mainexseq}) splits.
Let $1\in\wt{C}_0$ be a section projecting to a trivialization of $L$.
It is easy to see that the morphism $\OO_X[-1]\ra K$ corresponding
to the section $1^{r-1}\otimes \wt{d}(1)$ of $S^{r-1}\wt{C}_0\otimes C_1$
represents the map $\de$ from (\ref{symextri}). So from Lemma \ref{composlem}
we derive that $\wt{\tau}(1^{r-1}\otimes\wt{d}(1))=1$.
Together with the condition that the composition (\ref{zerocomp}) vanishes
this is equivalent to the equation
\begin{equation}
\lan\nu((x+\la\cdot 1)^{r-1}),\wt{d}(x+\la\cdot 1)\ran=\la^r
\end{equation} 
where $x\in C_0\sub\wt{C}_0$, the morphism
$\nu:S^{r-1}\wt{C}_0\ra C_1^{\vee}$
is induced by $\wt{\tau}$.
It follows that the section 
$$s_{\la}(x)=(\wt{d}(x+\la),\nu(x+\la\cdot 1))$$
of the orthogonal bundle $p^*C_1\oplus p^*C_1^{\vee}$
on $C_0$ satisfies $s_{\la}(x)\cdot s_{\la}(x)=\la^r$.
Applying Lemma \ref{mainlem} to the action of $s_{\la}$ on
the spinor bundle $\We^*p^*C_1^{\vee}$ we derive the vanishing of
localized Chern class corresponding to the isotropic section
$s_0$ obtained from $s_{\la}$ by setting $\la=0$. But the latter class is
precisely $c(\TT,b)$.
\ed

\section{Ramond factorization axiom}\label{Ramondsec}

Let the data $(\pi:\CC\ra X,\TT, b:\TT^r\ra\om_{\CC/X})$ 
be as in section \ref{vanishsec}.
Assume in addition that we have an $X$-point $\si:X\ra\CC$ which
is a nodal point of every fiber and that
$\wt{\pi}:\wt{\CC}\ra X$ is a fiberwise normalization of this point.
We denote by $n:\wt{\CC}\ra\CC$ the corresponding morphism and by
$\si_1,\si_2:X\ra\wt{\CC}$ two disjoint points that
project to $p\in\CC$. 
Finally, let us assume that $\TT$ is locally free at $\si$ and
that that map $b$ is an isomorphism at $\si$
(in \cite{JKV} this situation is referred to  as ``Ramond case'').

For every $\la\in k^*$ there is a natural line bundle
$\LL_{\la}$ on $\CC$ such that $n^*\LL_{\la}\simeq\OO_{\wt{\CC}}$
and the isomorphism $n^*\LL_{\la}|_{\si_1}\ra n^*\LL_{\la}|_{\si_2}$
corresponds to the multiplication by $\la$.
It is clear that $\LL_{\la\la'}\simeq\LL_{\la}\otimes\LL_{\la'}$.
In particular, if $\xi$ is an $r$-th root of unity then
$\LL_{\xi}^r\simeq\OO_{\CC}$. Therefore, we can twist the data 
$(\TT,b)$ by considering $\TT\otimes\LL_{\xi}$ and the
map $b_{\xi}:(\TT\otimes\LL_{\xi})^r\ra\om_{\CC/X}$ induced by
$b$ and the trivialization of $\LL_{\xi}^r$.

Now the Ramond case of Axiom 3 in \cite{JKV} is implied easily by
Theorem \ref{axiom1thm} together with the following result.

\begin{thm}\label{Ramondthm} One has
$$\sum_{\xi:\xi^r=1} c(\TT\otimes\LL_{\xi},b_{\xi})=0.$$
\end{thm}

Recall that the relative dualizing sheaves on $\CC/X$ and $\wt{\CC}/X$
are related by the isomorphism $n^*\om_{\CC/X}\simeq\om_{\wt{\CC}/X}(\si_1+\si_2)$
such that the following diagram is commutative
$$
\begin{diagram}
n^*\om_{\CC/X}|_{\si_1} &&\rTo && n^*\om_{\CC/X}|_{\si_2}\\
\dTo &&&&\dTo\\
\om_{\wt{\CC}/X}(\si_1+\si_2)|_{\si_1}&\rTo^{\Res}&\OO_X&
\lTo^{-\Res}&\om_{\wt{\CC}/X}(\si_1+\si_2)|_{\si_2}
\end{diagram}
$$
where the top arrow is the canonical isomorphism (the sign comes from the relation
$dx/x=-dy/y$ near the node $xy=0$). In particular, there is a canonical trivialization
of $\om_{\CC/X}|_{\si}$ such that
the boundary map $\de:\OO_X\simeq\om_{\CC/X}|_{\si}\ra R\pi_*\om_{\CC/X}[1]$ from the
exact triangle
$$R\pi_*\om_{\CC/X}\rTo R\pi_*n_*n^*\om_{\CC/X}\rTo\om_{\CC/X}|_{\si}\rTo^{\de}
R\pi_*\om_{\CC/X}[1]$$
satisfies
$\Tr\circ\de=\id$, where $\Tr:R\pi_*\om_{\CC/X}[1]\ra\OO_X$ is the trace map.

We can recover $R\pi_*\TT$ from $R\wt{\pi}_*\wt{\TT}$, where $\wt{\TT}=n^*\TT$,
together with the evaluation maps at $\si_1$ and $\si_2$. Namely, 
if we denote $L=\wt{\TT}|_{\si_1}\simeq\wt{\TT}|_{\si_2}$ then
there is an exact triangle
\begin{equation}\label{Ramondtriangle}
R\pi_*\TT\ra R\wt{\pi}_*\wt{\TT}\rTo^{\ev_{1}-\ev_{2}} L\ra R\pi_*\TT[1]
\end{equation}
where $\ev_{i}:R\wt{\pi}_*\wt{\TT}\ra L$ is the evaluation map at $\si_i$ ($i=1,2$).
Note that the morphism $b:\TT^r\ra\om_{\CC/X}$ induces a morphism
$$\wt{b}:\wt{\TT}^r\ra\om_{\wt{\CC}/X}(\si_1+\si_2).$$
Moreover, $\wt{b}$ is an isomorphism at $\si_1$ and $\si_2$, so restricting to either
of these points we get a trivialization of $L^r$ (the two trivializations are the same).
Passing to an \'etale cover of $X$ we can assume that $L$ itself is trivial.

Let $[C_0\stackrel{d}{\ra} C_1]$ be a complex of vector bundles representing $R\wt{\pi}_*\wt{\TT}$
with $C_1$ a direct sum of sufficiently negative powers of an ample line bundle on $X$. 
Then the evaluation maps 
$\ev_{1},\ev_{2}:R\wt{\pi}_*\wt{\TT}\ra L$ can be realized by morphisms 
$[C_0\ra C_1]\ra L$ in the homotopy category of complexes. Let 
$e_1,e_2:C_0\ra L$ be the corresponding morphisms (unique up to adding morphisms that factor
through $C_1$). Then we can choose a quasiisomorphism of
$R\pi_*\TT$ with the complex 
$$\Cone([C_0\ra C_1]\rTo^{e_1-e_2} L)[-1]=[C_0\rTo^{(d,e_1-e_2)}C_1\oplus L]$$
compatible with the triangle (\ref{Ramondtriangle}), where $\Cone(C\ra C')$ 
denotes the cone of a morphism of complexes $C\ra C'$. 

The triangle (\ref{Ramondtriangle}) is obtained by applying the functor $R\pi_*$ to
to the triangle 
$$\TT\ra n_*\wt{\TT}\rTo^{\ev_{1}-\ev_{2}} \si_*L\ra \TT[1]$$
on $\CC$. To understand the map $S^rR\pi_*\TT\ra R\pi_*\om_{\CC/X}$ we can use the
symmetric K\"unneth isomorphism $S^rR\pi_*\TT\simeq R\pi^{(r)}_*(\TT^{(r)})$, 
where $\TT^{(r)}$ is the $r$-th symmetric power of $\TT$ on $\CC^{(r)}$. 
The maps $\ev_1,\ev_2:n_*\wt{\TT}\ra\si_*L$ induce naturally the maps
$$\ev_1^r,\ev_2^r:(n_*\wt{\TT})^{(r)}\ra\si^r_*L,$$
where $\si^r:X\ra\CC^{(r)}$ is the $r$-tuple point of $\CC^{(r)}$ corresponding to $\si$.
Let us define a coherent sheaf on $\CC^{(r)}$ as follows: 
$$K:=\ker((n_*\wt{\TT})^{(r)}\rTo^{\ev_1^r-\ev_2^r}\si^r_*L^r).$$
Then we have a natural embedding $\TT^{(r)}\ra K$ which induces a morphism
$$\iota:S^rR\pi_*\TT\simeq R\pi^{(r)}_*\TT^{(r)}\ra R\pi^{(r)}_*K.$$
Let $\De:\CC\ra\CC^{(r)}$ be the diagonal embedding.
We claim that there is a natural morphism $K\ra\De_*\om_{\CC/X}$,
such that the composition of the induced map $\wt{\eta}:R\pi^{(r)}_*K\ra R\pi_*\om_{\CC/X}$
with $\iota$
coincides with the map $\eta:S^rR\pi_*\TT\ra R\pi_*\om_{\CC/X}$ induced by $b$. 
Indeed, $\De^*K$ maps to the kernel of the upper horizontal arrow in the commutative diagram
$$
\begin{diagram}
(n_*\wt{\TT})^{\otimes r}&\rTo^{\ev_1^r-\ev_2^r}&\si_*L^r\\
\dTo &&\dTo\\
n_*n^*\om_{\CC/X} &\rTo& \si_*(\om_{\CC/X}|_{\si})
\end{diagram}
$$
Therefore, we obtain the natural map from $\De^*K$ to the kernel of the lower
horizontal arrow in this diagram, i.e., a map $\De^*K\ra\om_{\CC/X}$. By adjunction
we get a morphism $K\ra\De_*\om_{\CC/X}$. The restriction of this map to the subsheaf
$\TT^{(r)}\sub K$ is the map induced by $b$ which implies our claim. We also
have a morphism of exact sequences
$$
\begin{diagram}
0 &\rTo& K&\rTo& (n_*\wt{\TT})^{(r)}&\rTo^{\ev_1^r-\ev_2^r}&\si^r_*L^r &\rTo& 0\\
&&\dTo&&\dTo &&\dTo\\
0 &\rTo&\De_*\om_{\CC/X}&\rTo&\De_*n_*n^*\om_{\CC/X} &\rTo&\si^r_*(\om_{\CC/X}|_{\si})&\rTo&0
\end{diagram}
$$
This implies the commutativity of the following diagram:
$$
\begin{diagram}
\OO_X\simeq &L^r &\rTo& R\pi^{(r)}_*K[1]\\
&\dTo&&\dTo^{\wt{\eta}[1]} \\
&\om_{\CC/X}|_{\si} &\rTo&R\pi_*\om_{\CC/X}[1]
\end{diagram}
$$
Therefore, the composition of the map 
$\wt{\tau}:=\Tr\circ\wt{\eta}[1]:R\pi^{(r)}_*K[1]\ra\OO_X$
with the natural map $\OO_X\simeq L^r\ra R\pi^{(r)}_*K[1]$ is equal to the identity.
On the other hand, since $\wt{\eta}\circ\iota=\eta$, it follows that the composition
$\wt{\tau}\circ\iota=\tau:S^rR\pi_*\TT[1]\ra\OO_X$ is exactly the map induced by $b$
(which is used in the definition of the class $c(\TT,b)$).

Note that $R\pi^{(r)}_*n_*\wt{\TT}\simeq S^rR\wt{\pi}_*\wt{\TT}$, so
the object $R\pi^{(r)}_*K$ fits into the distinguished triangle
$$R\pi^{(r)}_*K\ra S^r R\wt{\pi}_*\wt{\TT}\rTo^{\ev_1^r-\ev_2^r} L^r\ra R\pi^{(r)}_*K[1].$$
Therefore, it can be represented by the complex
$\Cone(S^r[C_0\ra C_1]\rTo^{e_1^r-e_2^r} L^r)[-1]$
in a way compatible with this triangle. Furthermore, the natural morphism
$S^rR\pi_*\TT\ra R\pi^{(r)}_*K$ is realized by the natural map of complexes
\begin{equation}\label{complexmap}
S^r[C_0\rTo^{(d,e_1-e_2)}C_1\oplus L]\ra \Cone(S^r[C_0\ra C_1]\rTo^{e_1^r-e_2^r} L^r)[-1]
\end{equation}
with the components $\id:S^r C_0\ra S^r C_0$, 
$$S^{r-1}C_0\otimes (C_1\oplus L)\ra (S^{r-1}C_0\otimes C_1)\oplus L^r:
x^{r-1}\otimes (y,z)\mapsto (x^{r-1}\otimes y, \sum_{i=0}^{r-1}e_1^i(x)e_2^{r-1-i}(x)z),$$
etc. Finally, we claim that for a suitable choice of the complex
$C_0\ra C_1$ (as in Proposition 4.7 of \cite{PV}) 
the map $\wt{\tau}:R\pi^{(r)}_*K[1]\ra\OO_X$ is
represented by the chain map of complexes $\Cone(S^r[C_0\ra C_1]\rTo^{e_1^r-e_2^r} L^r)\ra\OO_X$. 
This is a consequence of the following general result.

\begin{lem}\label{Ramondlem} 
Let $g:A\ra B$, $f:B\ra C$ be a pair of maps in the homotopy category $\KK$ of
some abelian category and let $\DD$ be the corresponding derived category.
Consider the subsets $H_{\KK}(f)\sub\Hom_{\KK}(\Cone(g),C)$ and $H_{\DD}(f)\sub\Hom_{\DD}(\Cone(g),C)$
consisting of morphisms $\Cone(g)\ra C$ such that their composition with the canonical
morphism $i:B\ra\Cone(g)$ is equal to $f$ (in $\KK$ and $\DD$ respectively).
Assume that the map
$$\Hom_{\KK}(A,C)\ra\Hom_{\DD}(A,C)$$
is injective and the map
$$\Hom_{\KK}(A[1],C)\ra\Hom_{\DD}(A[1],C)$$
is surjective. Then the natural map
$$\kappa:H_{\KK}(f)\ra H_{\DD}(f)$$
is surjective.
\end{lem}

\Pf . Let us denote by $\pi:\Cone(g)\ra A[1]$ the canonical chain map.
If the set $H_{\DD}(g)$ is empty then the assertion is clear, so we can assume
that $H_{\DD}(g)\neq\emptyset$. Then the composition $f\circ g:A\ra C$ becomes zero in
the derived category. By our assumption the natural map
$\Hom_{\KK}(A,C)\ra\Hom_{\DD}(A,C)$ is injective, hence $f\circ g$ is homotopic to zero.
Every homotopy $h$ from $g\circ f$ to $0$ induces naturally a chain map
$\Cone(h):\Cone(g)\ra C$ which coincides with $f$ on the subcomplex
$i(B)\sub\Cone(g)$. In fact, it is easy to see that the map $h\mapsto\Cone(h)$ is
a bijection between homotopies from $g\circ f$ to $0$ and chain maps
$\Cone(g)\ra C$ extending $f$ on $B$. If we have two homotopies $h_1,h_2$ from
$g\circ f$ to $0$ then the difference $h_1-h_2$ gives a chain map from $A[1]$ to $C$.
It is easy to see that 
$$\Cone(h_1)-\Cone(h_2)=(h_1-h_2)\circ\pi.$$

Now let $\ga\in H_{\DD}(g)$ be any element. Let us pick a homotopy $h_0$ from
$g\circ f$ to $0$. Then the homotopy class $[\Cone(h_0)]$ is an element of $H_{\KK}(f)$.
The composition of $\kappa([\Cone(h_0)])-\ga$ with $i$ vanishes in the derived category,
hence we have $\kappa([\Cone(h_0)])-\ga=\b\circ\pi$ for some
$\b\in\Hom_{\DD}(A[1],C)$. By our assumption there exists a chain map $\wt{\b}:A[1]\ra C$ 
representing $\b$. Then $h=h_0-\wt{\b}$ is another homotopy from $g\circ f$ to $0$.
We have 
$$\kappa([\Cone(h)])=\kappa([\Cone(h_0)-\wt{\b}\circ\pi])=\kappa([\Cone(h_0)])-\b\circ\pi=\ga.$$ 
\ed

We apply the above lemma to $A=S^r[C_0\ra C_1]$, $B=L^r$ and $C=\OO_X$, where
$f:L^r\ra\OO_X$ is the canonical isomorphism. To satisfy the assumptions
of the lemma we choose the complex $C_0\ra C_1$ representing $R\wt{\pi}_*\wt{\TT}$
with $C_1$ a direct sum of sufficiently negative
powers of an ample line bundle (one has to argue as in Proposition 4.7 of \cite{PV}).
Hence, the map $\wt{\tau}$ is represented by a morphism in the homotopic category
$\Cone(S^r[C_0\ra C_1]\rTo^{e_1^r-e_2^r} L^r)\ra\OO_X$ that we still denote by $\wt{\tau}$.
The restriction of $\wt{\tau}$ to the subcomplex $L^r$ is equal to the canonical isomorphism
$L^r\ra\OO_X$, while its composition with the map (\ref{complexmap}) is the morphism
$$\tau:S^r[C_0\rTo^{(d,e_1-e_2)}C_1\oplus L]\ra \OO_X[-1]$$
that should be used for the computation of $c(\TT,b)$.
It follows that the restriction of the corresponding morphism
$$\tau:S^{r-1}C_0\otimes (C_1\oplus L)\ra\OO_X$$
to $S^{r-1}C_0\otimes L$ has form
$$\tau(x^{r-1}\otimes y)=\sum_{i=0}^{r-1} e_{1}(x)^i e_{2}(x)^{r-1-i}y\in L^r\simeq\OO_X$$
where $x\in C_0$, $y\in L$.
Hence, the corresponding isotropic section of $p^*(C_1\oplus L\oplus C_1^{\vee}\oplus L^{-1})$
(where $p:C_0\ra X$ is the projection) has form
$$
s(x)=(d(x),(e_1-e_2)(x),\nu(x),\sum_{i=0}^{r-1} e_{1}(x)^i e_{2}(x)^{r-1-i}),
$$
where the last component belongs to $L^{r-1}\simeq L^{-1}$, $\nu$ is given by
some morphism $S^{r-1}C_0\ra C_1^{\vee}$.

To compute the class corresponding to the twisted data
$(\TT\otimes\LL_{\xi},b_{\xi})$ for some $r$-th root of unity $\xi$ we simply have to
replace the pair $(e_1,e_2)$ by $(e_1,\xi e_2)$. Note that
this will not affect the definition of $K$ and of the morphism $\wt{\tau}$.
Hence the corresponding isotropic section of $p^*(C_1\oplus L\oplus C_1^{\vee}\oplus L^{-1})$
will take form
\begin{equation}
s_{\xi}(x)=
(d(x),(e_1-\xi e_2)(x),\nu(x),\sum_{i=0}^{r-1} \xi^{r-1-i}e_{1}(x)^i e_{2}(x)^{r-1-i}),
\end{equation}
for some $\nu:S^{r-1}C_0\ra C_1^{\vee}$.

Now we can finish the proof of Theorem \ref{Ramondthm}.
For every $\xi$ the class $c(\TT\otimes\LL_{\xi},b_{\xi})$ is equal to 
$$\td(C_1\oplus L)\cdot\ch^{C_0}_X(\We^*p^*(C_1^{\vee}\oplus L^{-1}),s_{\xi}),$$
where $s_{\xi}\in p^*(C_1\oplus L\oplus C_1^{\vee}\oplus L^{-1}), s_{\xi})$
is the isotropic section $s_{\xi}$ constructed above.
Let us set $f_{\xi}=e_1-\xi e_2$. We consider $(f_{\xi})$ as a collection of
sections of $p^*L$ on $C_0$. We have an orthogonal decomposition
$$p^*(C_1\oplus L\oplus C_1^{\vee}\oplus L^{-1})\simeq p^*(C_1\oplus C_1^{\vee})\oplus
p^*(L\oplus L^{-1}),$$
so that the section $s_{\xi}$ has components
$s_0=(d,\nu)\in p^*(C_1\oplus C_1^{\vee})$ and
$(f_{\xi},\prod_{\xi'\neq\xi} f_{\xi'})$.
Recall that we can trivialize $L$, so the spinor
bundle $\We^*p^*(C_1^{\vee}\oplus L^{-1})$
can be identified with $\We^*p^*C_1^{\vee}\oplus \We^*p^*C_1^{\vee}[1]$. 
Under this identification
the action of sections $s_{\xi}$ will take form of differentials $(d_i)$
in Lemma \ref{mainlem2}, where
the odd endomorphism $d$ of $\We^*p^*C_1^{\vee}$ is given by the action of $s_0$.
Now the assertion of the theorem follows from Lemma \ref{mainlem2}.

\end{document}